\documentclass[12pt]{article}
\usepackage{amsmath,amssymb,mathrsfs}
\usepackage{hyperref}
\usepackage[PostScript=dvips,s=5ex]{diag}
\usepackage[latin1]{inputenc}
\diagramstyle[heads=curlyvee,s=5ex]

\begin{document}
\title{Convexity of Momentum Maps: A Topological Analysis}

\vspace{-21mm}
\author{Wolfgang Rump and Jenny Santoso\\ {\footnotesize\it Institute for
Algebra and Number Theory, University of Stuttgart}\\[-2mm]
{\footnotesize\it Pfaffenwaldring 57, D-70550 Stuttgart, Germany}\\[4mm]
{\small\it 
%
Dedicated to Alan D. Weinstein, Dennis P. Sullivan, and in memory
of}\\ 
{\small\it 
Jerrold E. Marsden}\\[4mm]} 
%
\date{}
\maketitle
\setcounter{footnote}{-1}
\renewcommand{\thefootnote}{\Roman{footnote}}
\footnotetext{\hspace*{-1.5mm} 2000 \it Mathematics Subject
Classification. \rm Primary: 52A01, 53C23, 54C10, 53D20. Secondary:
\rm 54E18. \it \hspace*{1mm} Key words and phrases. \rm Convexity
space, monotone-light factorization, geodesic manifold, momentum
map, Lokal-global-Prinzip.} \thispagestyle{empty}

\newtheorem{prop}{Proposition}
\newtheorem{thm}{Theorem}
\newtheorem{lem}{Lemma}
\newtheorem{Definition}{Definition}
\renewcommand{\labelenumi}{\rm (\alph{enumi})}
\newcommand{\hs}{\hspace{2mm}}
\newcommand{\vsp}{\vspace{4ex}}
\newcommand{\vspc}{\vspace{-1ex}}
\newcommand{\hra}{\hookrightarrow}
\newcommand{\tra}{\twoheadrightarrow}
\newcommand{\md}{\mbox{-\bf mod}}
\newcommand{\Mod}{\mbox{-\bf Mod}}
\newcommand{\Mdd}{\mbox{\bf Mod}}
\newcommand{\mdd}{\mbox{\bf mod}}
\newcommand{\latt}{\mbox{-\bf lat}}
\newcommand{\Proj}{\mbox{\bf Proj}}
\newcommand{\Inj}{\mbox{\bf Inj}}
\newcommand{\Ab}{\mbox{\bf Ab}}
\newcommand{\CM}{\mbox{-\bf CM}}
\newcommand{\Prj}{\mbox{-\bf Proj}}
\newcommand{\prj}{\mbox{-\bf proj}}
\newcommand{\ra}{\rightarrow}
\newcommand{\eps}{\varepsilon}
\renewcommand{\epsilon}{\varepsilon}
\renewcommand{\phi}{\varphi}
\renewcommand{\hom}{\mbox{Hom}}
\newcommand{\ex}{\mbox{Ext}}
\newcommand{\rad}{\mbox{Rad}}
\renewcommand{\Im}{\mbox{Im}}
\newcommand{\oti}{\otimes}
\newcommand{\sig}{\sigma}
\newcommand{\en}{\mbox{End}}
\newcommand{\Qq}{\mbox{Q}}
\newcommand{\lra}{\longrightarrow}
\newcommand{\lras}{\mbox{ $\longrightarrow\raisebox{1mm}{\hspace{-6.5mm}$\sim$}\hspace{3mm}$}}
\newcommand{\Eq}{\Leftrightarrow}
\newcommand{\Equ}{\Longleftrightarrow}
\newcommand{\Ra}{\Rightarrow}
\newcommand{\A}{\mathbb{A}}
\newcommand{\N}{\mathbb{N}}
\newcommand{\Z}{\mathbb{Z}}
\newcommand{\R}{\mathbb{R}}
\newcommand{\K}{\mathbb{K}}
\newcommand{\F}{\mathbb{F}}
\newcommand{\C}{\mathbb{C}}
\newcommand{\B}{\mathfrak{B}}
\newcommand{\rk}{\mbox{r}}
\newcommand{\p}{\mathfrak{p}}
\newcommand{\q}{\mathfrak{q}}
\renewcommand{\k}{\mathfrak{k}}
\newcommand{\AAA}{\mathfrak{A}}
\newcommand{\Aa}{{\cal A}}
\newcommand{\Dd}{{\cal D}}
\newcommand{\Bb}{{\cal B}}
\newcommand{\Ll}{{\cal L}}
\newcommand{\Hh}{{\cal H}}
\newcommand{\Nn}{{\cal N}}
\newcommand{\Cc}{{\cal C}}
\newcommand{\Mm}{{\cal M}}
\newcommand{\Tt}{{\cal T}}
\newcommand{\Ff}{{\cal F}}
\newcommand{\Rr}{{\cal R}}
\renewcommand{\AA}{\mathscr{A}}
\newcommand{\DD}{\mathscr{D}}
\newcommand{\BB}{\mathscr{B}}
\newcommand{\LL}{\mathscr{L}}
\newcommand{\HH}{\mathscr{H}}
\newcommand{\NN}{\mathscr{N}}
\newcommand{\EE}{\mathscr{E}}
\newcommand{\CC}{\mathscr{C}}
\newcommand{\UU}{\mathscr{U}}
\newcommand{\MM}{\mathscr{M}}
\newcommand{\PP}{\mathscr{P}}
\newcommand{\II}{\mathscr{I}}
\newcommand{\TT}{\mathscr{T}}
\newcommand{\XX}{\mathscr{X}}
\newcommand{\YY}{\mathscr{Y}}
\newcommand{\KK}{\mathscr{K}}
\newcommand{\FF}{\mathscr{F}}
\newcommand{\RR}{\mathscr{R}}
\newcommand{\PPP}{\mathfrak{P}}
\newcommand{\BBB}{\mathfrak{B}}
\newcommand{\DDD}{\mathfrak{D}}
\newcommand{\UUU}{\mathfrak{U}}
\newcommand{\FFF}{\mathfrak{F}}
\newcommand{\VVV}{\mathfrak{V}}
\newcommand{\WWW}{\mathfrak{W}}
\newcommand{\setm}{\smallsetminus}
\renewcommand{\le}{\leqslant}
\renewcommand{\ge}{\geqslant}
\newcommand{\rat}{\rightarrowtail}
\newcommand{\op}{^{\mbox{\scriptsize op}}}
\newcommand{\pf}{\it Proof.\hs\rm}
\newcommand{\bx}{\hspace*{\fill} $\square$}
\newcommand{\ltra}{\lra\!\!\!\!\!\!\:\ra}
\newcommand{\subs}{\subset}
\newcommand{\sups}{\supset}
\newcommand{\subsn}{\subsetneq}
\newcommand{\nsubs}{\not\subset}
\newcommand{\pt}{\makebox[0pt][r]{\bf .\hspace{.5mm}}}
\newcommand{\KS}{Krull-Schmidt }
\renewcommand{\theta}{\vartheta}
\newcommand{\dpt}{\makebox[0mm][l]{.}}
\newcommand{\dpc}{\makebox[0mm][l]{,}}
\newcommand{\MA}{\mbox{\textsf{M}$(\AA)$}}

\newarrow{To} ----{->}
\newarrow{to} ----{>}
\newarrow{in} {>}---{>}
\newarrow{inc} C---{>}
\newarrow{Inc} C---{->}
\newarrow{on} ----{>>}
\newarrow{On} ----{->>}
\newarrow{eq} =====
\newarrow{dash} {}{dash}{}{dash}{>}
\newarrow{dsh} {}{dash}{}{dash}{}

{\small\bf Abstract.} {\small The Local-to-Global-Principle used in
the proof of convexity theorems for momentum maps has been extracted
as a statement of pure topology enriched with a structure of
convexity. We extend this principle to not necessarily closed maps
$f\colon X\ra Y$ where the convexity structure of the target space
$Y$ need not be based on a metric. Using a new factorization of $f$,
convexity of the image is proved without local fiber connectedness,
and for arbitrary connected spaces $X$.}

\vspace{5mm}
\noindent {\Large\bf Introduction}

\vspace{5mm} Convexity for momentum maps was discovered
independently by Atiyah \cite{At} and Guillemin-Sternberg \cite{GS}
in the case of a Hamiltonian torus action 
on a compact symplectic\footnote[1]{\small For an introduction to
symplectic manifolds, their group actions, and momentum maps, see,
e.g., \cite{Wei2,Sil,MR,AM}}
manifold $X$. It was proved that the image of the momentum map $\mu$
is a convex polytope, namely, the convex hull of $\mu(X^T)$, where
$X^T$ denotes the set of fixed points under the action of the torus
$T$. In this case, $\mu$ is open onto its image, and the fibers of
$\mu$ are compact and connected. Two years later, in 1984, Kirwan
\cite{Kir} (see also \cite{GS1}) extended this result to the action
of a compact connected Lie group $G$. Here the image of $\mu\colon
X\ra\mbox{Lie}(G)^\ast$ has to be restricted to a closed Weyl
chamber in a Cartan subalgebra of $\mbox{Lie}(G)$, i.e. a
fundamental domain of $G$ with respect to its coadjoint action on
$\mbox{Lie}(G)^\ast$. Equivalently, this amounts to a composition of
the momentum map $\mu$ with the projection onto the quotient space
$Y:=\mbox{Lie}(G)^\ast/G$ modulo the coadjoint action of $G$. Up to
this time, convexity of $\mu$ was proved by means of Morse theory,
applied to the components of $\mu$. This works well as long as $\mu$
is defined on a compact manifold~$X$.

\textheight22.9cm 

In 1988, Condevaux, Dazord, and Molino \cite{CDM} reproved these
results in an entirely new fashion. They factor out the connected
components of the fibers of $\mu$ to get a monotone-light
factorization $\mu\colon X\ra\widetilde{X}\ra Y$ (see \cite{Mic}).
If $\mu$ is proper, i.e. closed and with quasi-compact fibers, the
metric of $Y$ can be lifted to $\widetilde{X}$. Then a shortest path
between two points of $\widetilde{X}$ maps to a straight line in
$Y$, which proves the convexity of $\mu(X)$. Based on this method,
Hilgert, Neeb, and Plank \cite{HNP} extended Kirwan's result to
non-compact connected manifolds $X$ under the assumption that $\mu$
is proper.

After this invention, the proof of convexity now splits into two
parts: A geometric part where certain local convexity data have to
be verified, and a topological part, a kind of
``Lokal-global-Prinzip'' \cite{HNP} which proves global convexity
\`a la Condevaux, Dazord, and Molino.

A further step was taken by Birtea, Ortega, and Ratiu
\cite{BOR1,BOR2} who consider a closed, not necessarily proper map
$\mu\colon X\ra\widetilde{X}\ra Y$, defined on a normal, first
countable, arcwise connected Hausdorff space $X$. The map $\mu$ has
to be locally open onto its image, locally fiber connected, having
local convexity data. Using Va\v{\i}n\v{s}te\v{\i}n's Lemma, they 
prove that the light part $\widetilde{X}\ra Y$ of $\mu$ is proper.
This yields global convexity of $\mu(X)$ for two almost disjoint
kinds of target spaces $Y$, either the dual of a Banach space
\cite{BOR2}
(which implies that the closed unit ball of $Y$ is weak$^\ast$ 
compact), or a complete locally compact length metric space $Y$
\cite{BOR1}. The second case applies to the cylinder-valued momentum
map \cite{OR1,OR2}, another invention of Condevaux, Dazord, and
Molino \cite{CDM}: For a symplectic manifold $(X,\omega)$, the
2-form $\omega$ gives rise to a flat connection on the trivial
principal fiber bundle $X\times\mbox{Lie}(G)^\ast$ with holonomy
group $H$. The cylinder-valued momentum map $\overline{\mu}$ is
obtained from $\mu$ by factoring out $\overline{H}$ from the target
space $Y$. The new target space $\overline{\mu}(X)=Y/\overline{H}$
is a cylinder, hence geodesics on it may differ from shortest paths.
The convexity theorem then states that $\overline{\mu}(X)$ is {\em 
weakly convex}, i.e. any two points of $\overline{\mu}(X)$ are
connected by a geodesic arc.

In the present paper, we analyse the topological part of convexity, that is,
the passage from local to global convexity. We show that the
Lokal-global-Prinzip, as developed thus far, admits a substantial improvement
in at least three respects.

Firstly, we replace the monotone-light factorization $f\colon X\ra
\widetilde{X}\ra Y$ that was used for a momentum map $f=\mu$ by a new
factorization
$$f\colon X\stackrel{q^f}{\ltra} X^f\stackrel{f^{\#}}{\lra} Y$$
of any continuous map $f\colon X\ra Y$ which is locally open onto its image.
In a sense, $X^f$ is closer
to $Y$ than the leaf space $\widetilde{X}$ since $q^f\colon X\ra X^f$
factors through the monotone part $X\ra\widetilde{X}$ of $f$.
We show that $q^f$ is an open surjection, while $X^f$
admits a basis of open sets $U$ such that $f^{\#}$ maps $U$ homeomorphically
onto a subspace of $Y$ (Proposition~\ref{p5}).
Therefore, $f^{\#}$ can take the r\^ole of the light part of $f$, which means
that we can drop the assumption that $f$ (the momentum map) is locally fiber
connected.

Secondly, we concentrate on the target space $Y$ instead of $X$ to
derive the desired properties of $X^f$. In this way, the various assumptions
on $X$ boil down to a single one, namely, its connectedness as a topological
space. Nevertheless, we need no extra assumptions on the target space $Y$.

Thirdly, we merely assume that the map $f^{\#}$ is closed, a much weaker
condition than the closedness of $f$. Even the light part of $f$
need not be closed. For example, $f^{\#}$ is trivial for a local homeomorphism
$f$ - a light map which need not be closed, and with fibers of arbitrary
size. Using the properties of $Y$, we prove that the fibers of $f^{\#}$ are
finite (Proposition~\ref{p10}), so that the convexity structure of $Y$ can be
lifted along $f^{\#}$ (Theorem~\ref{t2}).

To make the interaction between convexity and topology more visible,
we untie the Lokal-global-Prinzip from its metric context by means
of a general concept of convexity, which might be of interest in
itself. This also unifies the two above mentioned types of target
space considered in \cite{BOR1} and \cite{BOR2}. In the linear case
\cite{BOR2}, the target space $Y$ may be an arbitrary (not
necessarily complete) metrizable locally convex space instead of a
dual Banach space. (Metrizability can be weakened by the condition
that $Y$ does not contain a locally convex direct sum
$\R^{(\aleph_0)}$ as a subspace.) In general, geodesics in our
(non-linear) target space $Y$ are one-dimensional continua which
need not be metrizable.

In previous versions of the Lokal-global-Prinzip, geodesic arcs or
connecting lines between two points of the target space $Y$ are
obtained by a metric on $Y$. Without a concept of length, of course,
geodesics are no longer available by shortening of arcs in the
spirit of the Hopf-Rinow Theorem. Instead, we obtain geodesics by 
continued {\em straightening}, using a local convexity structure. In
other words, we deal with a ``manifold'', that is, a Hausdorff space
$Y$ covered by open subspaces $U$ with an additional structure of
convexity. The axioms of such a {\em convexity space} $U$ are very
simple: For any pair of points $x,y\in U$, there is a minimal
connected subset $C(x,y)$ containing $x$ and $y$, varying
continuously with the end points. In a topological vector space,
$C(x,y)$ is just the line segment between $x$ and $y$, while in a
uniquely geodesic space, $C(x,y)$ is the unique shortest path
between $x$ and $y$. With respect to the $C(x,y)$, there is a
natural concept of convexity, and for a convexity space $U$, we just
require that the $C(x,y)$ are convex and that $U$ has a basis of
convex open sets (see Definition~\ref{d1}).

If convexity is given by a metric, straightening and shortening of
arcs leads to the same result, namely, a geodesic of minimal length.
For a non-metrizable arc $A$ between two points $x$ and $y$, there
is a substitute for the length of $A$, namely, the closed convex
hull $\overline{C(A)}$ which is diminished by straightening. As a
first step, an inscribed ``line path'' $L$ (in a geodesic sense)
satisfies $\overline{C(L)}\subs\overline{C(A)}$, and
$\overline{C(L)}$ is the closed convex hull of the finitely many
extreme points of $L$. For a given line path $L$ between $x$ and
$y$, assume that the closed convex hull $\overline{C(L)}$ is
compact. Using Zorn's Lemma, we minimize the connected set 
$\overline{C(L)}$ to a compact convex set $C$ with $x,y\in C$. In
contrast to the Hopf-Rinow situation, where the shortening of $L$ is
achieved via the Arzel\`a-Ascoli Theorem, the straightening method 
needs the compactness of $\overline{C(L)}$ to show that
connectedness carries over to $C$. By the local convexity structure,
it then follows that $C$ contains a line path $L_0$ between $x$ and
$y$. Thus if $C=L_0$, the line path $L_0$ must be a geodesic.

So we require two properties to get the straightening process work:
First, the closed convex hull of a finite set must be compact; second, a
minimal compact connected convex set $C$ containing $x$ and $y$ has to be
a geodesic.

To establish a Lokal-global-Prinzip for continuous maps $X\ra Y$,
possible self-intersections 
of the arcs to be straightened have to
be taken into account. Precisely, this means that closed convex
subsets of $Y$ have to be replaced by {\em \'etale} maps, i.e.
closed locally convex maps $e\colon C\ra Y$, such that the connected
space $C$ admits a covering by open sets mapped homeomorphically
onto convex subsets of $Y$. We call $Y$ a {\em geodesic manifold} if
the above two properties hold with an adaption to \'etale maps
$e\colon C\ra Y$, that is, the second property now states that if
$C$ is compact and minimal with respect to $x,y\in C$, then $e$ can
be regarded as a geodesic with possible self-intersections. (Such a
geodesic is transversal if and only if $e=e^{\#}$.) If the charts
$U$ of $Y$ are regular Hausdorff spaces which satisfy a certain
finiteness condition (see Definition~\ref{d2}) which holds, for
example, if $U$ is either locally compact or first countable, we
call $Y$ a {\em geodesic $q$-manifold} (the ``q'' refers to the
finiteness condition). Obvious examples of geodesic $q$-manifolds
are complete locally compact length metric spaces, or metrizable
locally convex topological linear spaces (Examples~6 and 7). Our
main result consists in the following

\vspace{3mm}
\noindent {\bf Lokal-global-Prinzip.} \it Let $f\colon X\ra Y$ be a locally
convex continuous map from a connected topological space $X$ to a geodesic
$q$-manifold $Y$. Assume that $f^{\#}$ is closed. Then any two points of
$f(X)$ are connected by a geodesic arc.      \rm

For an inclusion map $f\colon C\hra Y$, the conditions on $f$ turn
into the assumptions of the Tietze-Nakajima Theorem (see  
\cite{Nak}), i.e. the subset $C$ is closed, connected, and locally
convex. Thus in case of a locally convex topological vector space
$Y$, the result for $C\hra Y$ yields Klee's Convexity Theorem 
\cite{Klee}, while for a complete Riemannian manifold $Y$, we get a
Theorem of Bangert \cite{Ban}. 

\section{Convexity spaces}

Let $X$ be a Hausdorff space. We endow the power set $\PPP(X)$ with a topology
as follows. For any open set $U$ of $X$, define
\begin{equation}\label{1}
\widetilde{U}:=\{ C\in\PPP(X)\: |\: C\subs U\}.
\end{equation}
The collection $\BBB$ of sets (\ref{1}) is closed under finite intersection.
We take $\BBB$ as a basis of open sets for the topology of $\PPP(X)$.
\begin{Definition}\pt\label{d1}
\rm Let $X$ be a Hausdorff space together with a continuous map
\begin{equation}\label{2}
C\colon X\times X\ra\PPP(X).
\end{equation}
We call a subset $A\subs X$ {\em convex} if $C(x,y)\subs A$ holds for all
$x,y\in A$. We say that $X$ is a {\em convexity space} with respect to a map
(\ref{2}) if the following are satisfied.

\setlength{\partopsep}{-3mm}\begin{enumerate}
\setlength{\parskip}{-1mm}
\item[(C1)] The $C(x,y)$ are convex for all $x,y\in X$.
\item[(C2)] The $C(x,y)$ are minimal among the connected sets $C\subs X$ with
$x,y\in C$.
\item[(C3)] $X$ has a basis of convex open sets.
\end{enumerate}
\end{Definition}

Note that (C1) implies that $C(y,x)\subs C(x,y)$. Hence $C$ is symmetric:
\begin{equation}\label{3}
C(x,y)=C(y,x).
\end{equation}
From (C2) we infer that
\begin{equation}\label{4}
C(x,x)=\{ x\}.
\end{equation}
Moreover, (C2) implies that every convexity space $X$ is connected. The
restriction of the map (\ref{2}) to a convex subset $A\subs X$ makes $A$ into
a convexity space. Hence (C3) implies that $X$ is locally connected.
\begin{lem}\pt\label{l1}
Let $X$ be a convexity space. For $x,y\in X$, the set $C(x,y)\setm\{ y\}$ is
connected.
\end{lem}

\pf Let $A$ be the connected component of $x$ in $C(x,y)\setm\{ y\}$. Since
$\{ y\}$ is closed, every $z\in C(x,y)\setm\{ y\}$ admits a convex
neighbourhood $U$
with $y\notin U$. Hence $C(x,y)\setm\{ y\}$ is locally connected, and thus $A$
is open in $C(x,y)$. Since $C(x,y)$ is connected, it follows that $A$ cannot be
closed in $C(x,y)$. Thus $y\in\overline{A}$, which shows that $A\cup\{ y\}$
is connected. By (C2), this gives $A\cup\{ y\}=C(x,y)$, whence $A=C(x,y)\setm
\{ y\}$. \bx

As a consequence, the $C(x,y)$ can be equipped with a natural ordering.
\begin{prop}\pt\label{p1}
Let $X$ be a convexity space. For $x,y\in X$, the set $C(x,y)$ is linearly
ordered by
\begin{equation}\label{5}
z\le t\; :\Equ\; z\in C(x,t)\;\Equ\; t\in C(z,y)
\end{equation}
for $z,t\in C(x,y)$.
\end{prop}

\pf For any $z\in C(x,y)$, the set $C(x,z)\cup C(z,y)$ is connected. Therefore,
(C1) and (C2) give
\begin{equation}\label{6}
C(x,y)=C(x,z)\cup C(z,y).
\end{equation}
To verify the second equivalence in (\ref{5}), it suffices to show that
$$z\in C(x,t)\;\Ra\; t\in C(z,y)$$
holds for $z,t\in C(x,y)$. By Eq. (\ref{6}), it is enough to prove the
implication
\begin{equation}\label{7}
z\in C(x,t)\setm\{ t\}\;\Ra\; t\notin C(x,z).
\end{equation}
Assume that $z\in C(x,t)\setm\{ t\}$. Then Eq. (\ref{4}) gives $x\in C(x,t)
\setm\{ t\}$. Hence Lemma~\ref{l1}
and (C2) yield $C(x,z)\subs C(x,t)\setm\{ t\}$, which proves (\ref{7}).
Clearly, the relation (\ref{5}) is reflexive and transitive. By (\ref{7}), it
is a partial order. Furthermore, (\ref{5}) and (\ref{6}) imply that it is a
linear order. \bx

Note that the ordering of $C(x,y)$ depends on the pair $(x,y)$ which determines
the initial choice $x\le y$. Thus as an ordered set, $C(y,x)$ is dual to
$C(x,y)$.

\noindent {\bf Example 1.} Let $\Omega$ be a linearly ordered set. A
subset $I$ of $\Omega$ is said to be an {\em interval} if $a\le c\le
b$ with $a,b\in I$ implies that $c\in I$. The intervals $\{
c\in\Omega\:|\:c<b\}$ and $\{ c\in\Omega\:|\:c>a\}$ with
$a,b\in\Omega$ form a sub-basis for the {\em order topology} of
$\Omega$. Note that an open set of $\Omega$ is a disjoint union of
open intervals. Therefore, $\Omega$ is connected if and only if it
is a {\em linear continuum}, i.e. if every partition $\Omega=I\sqcup
J$ into non-empty intervals $I,J$ determines a unique element
between $I$ and $J$. With the order topology, a linear continuum
$\Omega$ is a locally compact convexity space with
\begin{equation}\label{8}
C(x,y)=\{ z\in\Omega\: |\: x\le z\le y\}
\end{equation}
in case that $x\le y$. Here the convex sets of $\Omega$ are just the connected
sets of $\Omega$.

\noindent {\bf Example 2.} More generally, we define a {\em tree
continuum} to be a Hausdorff space $X$ for which every two points
$x,y\in X$ are contained in a smallest connected set $C(x,y)$ such
that each $C(x,y)$ is a linear continuum, and $X$ carries the finest
topology for which the inclusions $C(x,y)\hra X$ are continuous.
Thus $U\subs X$ is open if and only if every $x\in U$ is an
``algebraically inner'' point (see \cite{Koe}, \S 16.2), i.e. if for
each $y\in X\setm\{ x\}$, there exists some $z\in C(x,y)\setm\{x\}$
with $C(x,z)\setm\{z\}\subs U$. Then $X$ is a convexity space. For
example, every one-dimensional CW-complex without cycles is of this
type.

\noindent {\bf Example 3.} In the Euclidean plane $\R^2$, consider the
solution curves $c\colon\R\ra\R^2$
of the differential equation $y'=3y^{\frac{3}{2}}$ (including the
singular solution $c\colon x\mapsto\tbinom{x}{0}$). With the finest topology
making the solution curves continuous, $\R^2$ becomes a tree continuum.
Here every point of the singular line is a branching point of order 4.

The following lemma is well-known (see \cite{Wil}, Theorem~26.15). 
\begin{lem}\pt\label{l2}
Let $X$ be a connected topological space with an open covering $\UUU$. For any
pair of points $x,y\in X$, there is a finite sequence $U_1,\ldots,U_n\in\UUU$
with $x\in U_1$, $y\in U_n$, and $U_i\cap U_{i+1}\not=\varnothing$ for $i<n$.
\end{lem}

\begin{prop}\pt\label{p2}
Let $X$ be a convexity space. For $x,y\in X$, the subspace $C(x,y)$ is
compact and carries the order topology.
\end{prop}

\pf Let $C(x,y)=\bigcup\UUU$ be a covering by convex open sets. By
Lemma~\ref{l2}, there is a finite sequence $U_1,\ldots,U_n\in\UUU$
with $x\in U_1, y\in U_n$, and $U_i\cap U_{i+1}\not=\varnothing$ for $i<n$.
Hence $C(x,y)=U_1\cup\cdots\cup U_n$, which shows that $C(x,y)$ is compact.

\vspc
For $u<v$ in $C(x,y)$, the sets $C(x,u)$ and $C(v,y)$ are compact, hence
closed in $C(x,y)$. So the open intervals of $C(x,y)$ are open sets in
$C(x,y)$. Conversely, a convex open set in $C(x,y)$ is an interval
which must be an open interval since $C(x,y)$ is connected. \bx

Up to here, we have not used the continuity of the map (\ref{2}) in
Definition~\ref{d1}.
\begin{prop}\pt\label{p3}
Let $X$ be a convexity space. The closure of any convex set $A\subs X$ is
convex.
\end{prop}

\pf Let $A\subs X$ be a convex set, and let $x,y\in\overline{A}$ be
given. For any $z\in C(x,y)$, we have to show that
$z\in\overline{A}$. Suppose that there is a convex neighbourhood $W$
of $z$ with $W\cap A=\varnothing$. Then $z\not=x,y$. By
Proposition~\ref{p2}, there exist $u,v\in W\cap C(x,y)$ with
$u<z<v$. Since $C(x,u)$ and $C(v,y)$ are compact, there are disjoint
open sets $U,V$ in $X$ with $C(x,u)\subs U$ and $C(v,y)\subs V$
(see, e.g., \cite{Kel}, chap.~V, Theorem~8). Hence $C(x,y)\subs
U\cup V\cup W$. So there are neighbourhoods $U'\subs U$ of $x$ and
$V'\subs V$ of $y$ with $C(x',y')\subs U\cup V\cup W$ for all $x'\in
U'$ and $y'\in V'$. Choose $x',y' \in A$. Then $C(x',y')\subs A$,
which yields $C(x',y')\subs U\cup V$, where $x'\in U'\subs U$ and
$y'\in V'\subs V$, contrary to the connectedness of $C(x',y')$. \bx

\begin{Definition}\pt\label{d2}
\rm Let $X$ be a convexity space. Define a {\em star} in $X$ with {\em center}
$x\in X$ and {\em end set} $E\subs X\setm\{ x\}$ to be a subspace $S(x,E):=
\bigcup\{ C(x,z)\:|\:z\in E\}$ with $C(x,z)\cap C(x,z')=\{x\}$ for
different $z,z'\in E$ such that
$S(x,E)$ carries the finest topology which makes the embeddings $C(x,z)\hra
S(x,E)$ continuous for all $z\in E$. We call $X$ {\em star-finite} if every
closed star in $X$ has a finite end set.
\end{Definition}
Thus every star is a tree continuum (Example~2). Recall that a topological
space $X$ is said to be a {\em $q$-space}
\cite{Mic1} if every point of $X$ has a sequence $(U_n)_{n\in\N}$ of
neighbourhoods such that every sequence $(x_n)_{n\in\N}$ with $x_n\in U_n$
admits an accumulation point. For example, every locally compact space, and
every first countable space $X$ is a $q$-space.

\begin{prop}\pt\label{p4}
Let $X$ be a convexity space which is a $q$-space. Then $X$ is star-finite.
\end{prop}

\pf Let $S(x,E)$ be a closed star in $X$, and let $(U_n)_{n\in\N}$ be
a sequence of neighbourhoods of $x$ such that every sequence
$(x_n)_{n\in\N}$ with $x_n\in U_n$ has an accumulation point. Suppose that $E$
is infinite. Since $U_n\cap C(x,z)\not=\{x\}$ for all $n\in\N$ and $z\in E$,
we find a subset $\{z_n\: |\:n\in\N\}$ of $E$ and a sequence $(x_n)_{n\in\N}$
with $x\not=x_n\in C(x,z_n)\cap U_n$. Thus $(x_n)_{n\in\N}$ has an
accumulation point $z$. Because of the star-topology, $z$ cannot belong to
$S(x,E)$, contrary to the assumption that $S(x,E)$ is closed. \bx

\noindent {\bf Example 4.} A topological vector space $X$ is a convexity
space with respect to straight line segments if and only if $X$ is locally
convex. Moreover, a locally convex space $X$ is star-finite if and only if $X$
does not contain a locally convex direct sum $\R^{(\aleph_0)}$ as a subspace.
In fact, every subspace $\bigoplus_{x\in E}\R x$ of $X$ with
$|E|=\aleph_0$ is complete (\cite{Sch}, II.6.2) and gives rise to a
closed star $S(0,E)$. Conversely, let $S(x,E)\subs X$ be a closed star with
$E$ infinite. Since finite dimensional subspaces of $X$ are star-finite by
Proposition~\ref{p4}, we can assume that $E$ is linearly independent and
$x=0$. Then the subspace $\bigoplus_{x\in E}\R x$ of $X$ is a locally convex
direct sum.

\vspc
Note that every metrizable locally convex space  $X$ is first countable
(\cite{Sch}, I, Theorem~6.1), hence star-finite by Proposition~\ref{p4}.

\section{Local openness onto the image}

For a topological space $X$, the infinitesimal structure at a point $x$ is
given by the set $\DDD_x$ of filters on $X$ which converge to $x$. Let
$\FFF(X)$ denote the set of all filters on $X$. We make $\FFF(X)$ into a
topological space with a basis of open sets
\begin{equation}\label{9}
\widetilde{U}:=\{\alpha\in\FFF(X)\: |\: U\in\alpha\},
\end{equation}
where $U$ runs through the class of open sets in $X$. Every continuous map
$f\colon X\ra Y$ induces a map $\FFF(f)\colon\FFF(X)\ra\FFF(Y)$. For an open
set $V$ in $Y$, we have
\begin{equation}\label{10}
\FFF(f)^{-1}(\widetilde{V})=\widetilde{f^{-1}(V)},
\end{equation}
which shows that $\FFF(f)$ is continuous.
Consider the subspace
\begin{equation}\label{11}
\DDD(X):=\{(x,\alpha)\in X\times\FFF(X)\: |\:\alpha\in\DDD_x\}
\end{equation}
of $X\times\FFF(X)$.
Note that for every $x\in X$, the neighbourhood filter $\UU(x)$ of $x$ is the
coarsest filter in $\DDD_x$. Thus, regarding $\DDD_x$ as a subset of $\DDD(X)$,
we get a pair of continuous maps
\begin{equation}\label{12}
\begin{diagram}[midshaft]
X & \rTo^{\mbox{\scriptsize $\UU$}} & \DDD(X) & \rOn^{\mbox{\scriptsize lim}}
& X
\end{diagram}
\end{equation}
with $\mbox{lim}(x,\alpha):=x$ and $\mbox{lim}\circ\UU=1_X$. In particular,
$\DDD_x=\mbox{lim}^{-1}(x)$.

For a continuous map $f\colon X\ra Y$, the local behaviour at
$x\in X$ is given by the induced map $\DDD_x f\colon\DDD_x\ra
\DDD_{f(x)}$. Thus we get an endofunctor $\DDD\colon\mbox{\bf Top}\ra
\mbox{\bf Top}$ of the category \mbox{\bf Top} of topological spaces with
continuous maps as morphisms. The functor $\DDD$ is augmented by the natural
transformation $\mbox{lim}\colon\DDD\ra 1$. On the other hand, the equation
$\UU\circ f=\DDD(f)\circ\UU$ holds if and only if $f$ is open.
\begin{Definition}\pt\label{d3}
\rm A continuous map $f\colon X\ra Y$ between topological spaces is said to
be {\em locally open onto its image} \cite{Ben} if every $x\in X$ admits an
open neighbourhood $U$ such that the induced map $U\tra f(U)$ is open onto the
subspace $f(U)$ of $Y$. We call $f$ {\em filtered} if $f$ is locally open onto
its image and $\DDD(f)\circ\UU$ is injective.
\end{Definition}
For example, the identity map $1_X\colon X\ra X$ is filtered if and
only if every point of $X$ is determined by its neighbourhood
filter, i.e. if $X$ is a $T_0$-space. The following structure
theorem holds for continuous maps which are locally open onto their
image.
\begin{prop}\pt\label{p5}
Let $f\colon X\ra Y$ be a continuous map which is locally open onto its image.
Up to isomorphism, there is a unique factorization $f=pq$ in $\mbox{\bf Top}$
into an open surjection $q$ and a
filtered map $p$. If $f$ is filtered, then every point $x\in X$ has an open
neighbourhood which is mapped homeomorphically onto a subspace of $Y$.
\end{prop}

\pf Consider the following commutative diagram
$$\begin{diagram}[s=6ex,midshaft]
1\colon\hspace*{-10mm} & X & \rTo^{\mbox{$\UU$}} & \DDD(X) & \rTo^{\mbox{lim}}
& X\\
& \dOn>{\mbox{$q^f$}} & & \dTo>{\mbox{$\DDD(f)$}} & & \dTo>{\mbox{$f$}}\\
f^{\#}\colon\hspace*{-6mm} & X^f & \rInc^{\mbox{$e$}} & \DDD(Y) &
\rTo^{\mbox{lim}} & Y\dpc\\
\end{diagram} $$
where $X^f$ is the image of $\DDD(f)\circ\UU$, regarded as a
quotient space of $X$, and $f^{\#}:=\mbox{lim}\circ e$. We will
prove that $f=f^{\#}\circ q^f$ gives the desired factorization. Let
us show first that $q^f$ is open. Thus let $U$ be an open set of
$X$. We have to verify that $(q^f)^{-1}q^f(U)$ is open in $X$. Since
$f$ is locally open onto its image, we can assume that the induced
map $U\tra f(U)$ is open. Let $x\in (q^f)^{-1}q^f(U)$ be given. Then
$q^f(x)\in q^f(U)$. So there exists some $y\in U$ with
$q^f(x)=q^f(y)$, i.e. $f(x)=f(y)$ and $f(\UU(x))=f(\UU(y))$. Hence
there is an open neighbourhood $V\in\UU(x)$ with $f(V)\subs f(U)$.
Again, we can assume that the induced map $V\tra f(V)$ is open.
Furthermore, there is an open neighbourhood $U'\subs U$ of $y$ with
$f(U')\subs f(V)$, and $f(U')$ is open in $f(U)$, hence in $f(V)$.
Therefore, $V':=V\cap f^{-1}(f(U'))$ is an open neighbourhood of $x$
with $f(V')=f(U')$.

\vspc
For any $x'\in V'$, there is a point $y'\in U'$ with $f(x')=f(y')$. So the
continuity of $f$ implies that $f(\UU(x'))=f(\UU(y'))$, which gives
$q^f(x')=q^f(y')$, and thus $V'\subs (q^f)^{-1}q^f(U')\subs (q^f)^{-1}q^f(U)$.
This proves that
$q^f$ is open. Consequently, $f^{\#}$ is locally open onto its image.

\vspc
Since $q^f$ is open, we have a commutative diagram
$$\begin{diagram}[w=8ex,h=5ex,midshaft]
X & \rOn^{\mbox{$q^f$}} & X^f\\
\dTo>{\mbox{$\UU$}} & & \dTo>{\mbox{$\UU$}}\\
\DDD(X) & \rTo^{\mbox{$\DDD(q^f)$}} & \DDD(X^f).\\
\end{diagram} $$
Hence $\DDD(f^{\#})\circ\UU\circ q^f=\DDD(f^{\#})\circ\DDD(q^f)
\circ\UU=\DDD(f)\circ\UU=e\circ q^f$. Therefore, $\DDD(f^{\#})\circ\UU=e$,
which implies that $f^{\#}$ is filtered.

\vspc
Now let $f=pq=p'q'$ be two factorizations with $p,p'$ filtered and $q,q'$
open. Then $\DDD(p')\circ\UU\circ q'=\DDD(p')\circ\DDD(q')\circ\UU=\DDD(p)
\circ\DDD(q)\circ\UU=\DDD(p)\circ\UU\circ q$. Since $\DDD(p')\circ\UU$ is
injective, there exists a map $e\colon E\ra E'$ with $q'=eq$. Since $q$ is
open, the map $e$ is continuous. So we get a commutative diagram
$$\begin{diagram}[s=5ex]
X & \rOn^{\mbox{$q$}} & E & \rTo^{\mbox{$p$}} & Y\\
\deq & & \dTo>{\mbox{$e$}} & & \deq\\
X & \rOn^{\mbox{$q'$}} & E' & \rTo^{\mbox{$p'$}} & Y\\
\end{diagram} $$
in \mbox{\bf Top}. By symmetry, we find a continuous map $e'\colon E'\ra E$
with $q=e'q'$ and
$p'=pe'$. Since $q$ and $q'$ are surjective, $e$ must be a homeomorphism. This
proves the uniqueness of the factorization.

\vspc
Finally, let $f\colon X\ra Y$ be filtered. For a given $x\in X$, let $U$ be an
open neighbourhood such that the induced map $r\colon U\tra f(U)$ is open.
Since $i\colon U\hra X$ is open, we have a commutative diagram
$$\begin{diagram}[w=8ex]
& & X & \rTo^{\mbox{$\UU$}} & \DDD(X) & & \\
& \ldTo \makebox[0pt][l]{\hspace*{-2.5mm}\raisebox{2mm}{$f$}} &
\uInc>{\mbox{$i$}} & & \uTo>{\mbox{$\DDD(i)$}} & \rdTo
\makebox[0pt][l]{\hspace*{-.5mm}\raisebox{1.5mm}{$\DDD(f)$}} & \\
\mbox{\hspace*{4.9mm} $Y$\rule[-.77mm]{0pt}{4.1mm}} & & U &
\rTo^{\mbox{$\UU$}} & \DDD(U) & & \DDD(Y) \\
& \luInc \makebox[0pt][l]{\hspace*{-2.5mm}\raisebox{-2mm}{$j$}} &
\dOn>{\mbox{$r$}} & & \dTo>{\mbox{$\DDD(r)$}} & \ruTo
\makebox[0pt][l]{\hspace*{-.5mm}\raisebox{-2mm}{$\DDD(j)$}} & \\
& & f(U) & \rTo^{\mbox{$\UU$}} & \DDD(f(U)) & & \\
\end{diagram} $$
which shows that $\DDD(j)\circ\UU\circ r=\DDD(f)\circ\UU\circ i$ is injective.
Hence $r$ is injective. \bx

In the sequel, we keep the notation of Proposition~\ref{p5} and write
\begin{equation}\label{13}
f\colon X\stackrel{q^f}{\ltra} X^f\stackrel{f^{\#}}{\lra} Y
\end{equation}
for the factorization of a map $f$ which is locally open onto its image.

\noindent {\bf Remarks. 1.} Although the factorization (\ref{13}) is
unique up to isomorphism, it does not give rise to a factorization
system \cite{FK,CHK}, i.e. a pair $(\EE,\MM)$ of subcategories such
that every commutative square
\begin{equation}\label{14}
\begin{diagram}
E_1 & \rTo^{\mbox{$f_1$}} & M_1\\
\dTo<{\mbox{$e$}} & \makebox[0pt][l]{\hspace*{-5.5mm} $d$}\rudash &
\dTo>{\mbox{$m$}}\\
E_0 & \rTo^{\mbox{$f_0$}} & M_0\\
\end{diagram}
\end{equation}
with $e\in\EE$ and $m\in\MM$ admits a unique diagonal $d$ with
$f_1=de$ and $f_0=md$ (see \cite{Her}, Proposition~1.4). Apart from
the fact that local openness onto the image is not closed under
composition (consider the maps $\R\stackrel{i}{\hra}\R^2
\stackrel{p}{\tra}\R$ with $i(x)=\binom{x}{x^3-3x}$ and
$p\colon\binom{x}{y} \mapsto y$), there cannot be a factorization
system since open surjections are not stable under pushout (take,
e.g., the pushout of the open surjection $\R\tra\{0\}$ and the
inclusion $\R\hra\R^2$).

{\bf 2.} If $f\colon X\ra Y$ is locally open onto its image and
locally fiber connected \cite{Ben,HNP}, the Lemma of Benoist 
(\cite{Ben}, Lemma~3.7) states that the monotone part $\pi$ of the
monotone-light factorization $f=\widetilde{f}\circ\pi$ is open. Here
the local fiber-connectedness of $f$ implies that $\pi$ is locally
open onto its image. Hence $\pi=q^\pi$ is open by
Proposition~\ref{p5}. In general, $q^f$ always factors through
$\pi$, but the two factorizations need not be isomorphic. For
example, a local homeomorphism $f\colon X\tra Y$ is open, but its
fibers are discrete.

\section{Convexity of maps}

In this brief section, we introduce local convexity and extend this concept
from subsets to continuous maps (cf. \cite{KB} for a notion of convex maps
in terms of paths).
\begin{Definition}\pt\label{d4}
\rm Let $X$ be a topological space. We define a {\em local convexity
structure}
on $X$ to be an open covering $X=\bigcup\UUU$ by convexity spaces $U\in\UUU$
(with the induced topology) such that for any $U\in\UUU$, every convex open
subspace of $U$ belongs to $\UUU$ (as a convexity space). We call a subset
$C\subs X$ {\em convex} if $C\cap U$ is convex for all $U\in\UUU$. We say
that $C$ is {\em locally convex} if every $z\in C$ admits a neighbourhood
$U\in\UUU$ such that $C\cap U$ is convex.
\end{Definition}
The covering $\UUU$ will be referred to as the {\em atlas} of the
local convexity structure. In the special case $X\in\UUU$, the atlas $\UUU$
just consists of the convex open sets of a convexity space $X$.

In contrast to local convexity, our concept of convexity refers to all
sets in $\UUU$. So the
intersection of convex sets is convex, and every subset $A\subs X$ admits a
{\em convex hull} $C(A)$, that is, a smallest convex set $C\sups A$. The next
proposition generalizes Proposition~\ref{p3}.
\begin{prop}\pt\label{p6}
Let $X$ be a topological space with a local convexity structure $\UUU$. The
closure of any convex set $A\subs X$ is convex.
\end{prop}

\pf For every $U\in\UUU$, we have $\overline{A}\cap U=\overline{A\cap U}
\cap U$. This set is convex by Proposition~\ref{p3}. Hence $\overline{A}$ is
convex. \bx

Definition~\ref{d4} admits a natural extension to continuous maps.

\begin{Definition}\pt\label{d5}
\rm Let $f\colon X\ra Y$ be a continuous map between topological spaces, where
$Y$ has a local convexity structure $\VVV$. We call $f$
{\em locally convex} if every $x\in X$ admits an open neighbourhood $U$ such
that the induced map $U\tra f(U)$ is open, and $f(U)$ is a convex subspace
of some $V\in\VVV$.
\end{Definition}

\noindent {\bf Remarks. 1.} A subset $A\subs Y$ is locally convex if and only
if the inclusion map $A\hra Y$ is locally convex.

\vspc
{\bf 2.} The open neighbourhood $U$ of $x$ in Definition~\ref{d5} can be chosen
arbitrarily small. In fact, let $U'\subs U$ be any smaller open neighbourhood
of $x$. Then $f(U')$ is an open subset of $f(U)$. Hence there exists some
$V'\in\VVV$ with $f(x)\in V'\cap f(U)\subs f(U')$. Thus $U'':=
U'\cap f^{-1}(V')$ is an open neighbourhood of $x$ with $f(U'')=V'\cap f(U')=
V'\cap f(U)$, which is a convex subspace of $V'$.

\vspc {\bf 3.} If $X$ is a connected Hausdorff space and $Y$ a
length metric space \cite{BrH,Gr}, a continuous map $f\colon X\ra Y$
is locally convex if and only if $f$ is locally open onto its image
and has local convexity data in the sense of \cite{BOR1}.

\begin{prop}\pt\label{p7}
Let $f\colon X\ra Y$ be a continuous map between topological spaces, where $Y$
has a local convexity structure $\VVV$. If $f$ is
locally convex, then $f^{\#}$ is locally convex.
\end{prop}

\pf Assume that $f$ is locally convex, and let $U$ be an open neighbourhood of
$x\in X$ such that the induced map $U\tra f(U)$ is open onto a convex
subspace of some $V\in\VVV$. Since $q^f$ is open by
Proposition~\ref{p5}, this property of $U$ carries over to the neighbourhood
$q^f(U)$ of $q^f(x)$. Hence $f^{\#}$ is locally convex. \bx

\section{Geodesic manifolds}

In this section, we introduce a general concept of geodesic which does not
refer to any kind of metric.

\begin{Definition}\pt\label{d6}
\rm Let $Y$ be a topological space with a local convexity structure $\VVV$,
and let $e\colon C\ra Y$ be a continuous map with a connected topological space
$C$. By $\VVV_e$ we denote the set of all open sets $U$ in $C$ which are mapped
homeomorphically onto a convex subspace of some $V\in\VVV$. We call $e$
{\em \'etale} if $e$ is closed and $\VVV_e$ covers $C$. We say that
$e\colon C\ra Y$ is {\em generated} by a subset $F\subs C$ if
there is no closed connected subspace $A\subsn C$ with $F\subs A$ such that
$e(U\cap A)$ is convex for all $U\in\VVV_e$.
\end{Definition}
In particular, \'etale maps are locally convex. Furthermore, every \'etale
map $e\colon C\ra Y$ induces a local convexity structure $\VVV_e$ on $C$. So
the condition (Definition~\ref{d6}) that $e(U\cap A)$ is convex for all
$U\in\VVV_e$ just states that $A$ is convex with respect to $\VVV_e$. If
$F\subs C$ is connected, then $C(F)$ is connected. Therefore, an
\'etale map $e\colon C\ra Y$ is generated by a connected set $F$ if and only
if $\overline{C(F)}=C$. Note that the composition of \'etale maps is \'etale.

\begin{Definition}\pt\label{d7}
\rm Let $Y$ be a Hausdorff space with a local convexity structure $\VVV$. We
call $Y$ a {\em geodesic manifold} if the following are satisfied.

\setlength{\partopsep}{-3mm}\begin{enumerate}
\setlength{\parskip}{-1mm}
\item[(G1)] For a finite set $F\subs Y$, the closure of $C(F)$ is compact.
\item[(G2)] If an \'etale map $e\colon C\ra Y$ with $C$ compact is generated by
$\{x,y\}\subs C$, then every connected set $A\subs C$ with $x,y\in A$
coincides with $C$.
\end{enumerate}
If, in addition, every $V\in\VVV$ is star-finite and regular (as a topological
space), we call $Y$ a {\em geodesic $q$-manifold}.
\end{Definition}
The letter ``$q$'' is reminiscent of Proposition~\ref{p4}. Since a geodesic
manifold $Y$ is locally connected, \cite{Bou}, chap.~I, 11.6,
Proposition~11, implies that $Y$ is the topological sum of its connected
components.

\begin{Definition}\pt\label{d8}
\rm Let $Y$ be a geodesic manifold. We define a {\em geodesic} in $Y$ to be
an \'etale map $e\colon C\ra Y$, generated by $\{ x,y\}\subs C$, where $C$ is
compact. The points
$e(x)$ and $e(y)$ will be called the {\em end points} of the geodesic.
\end{Definition}

More generally, we define a {\em line path} in $Y$ to be a continuous map
\mbox{$e\colon L\ra Y$}, where $L$ is a linear continuum (Example~1) with end
points $x_0$ and $x_n$ and a sequence of intermediate points $x_0<x_1<\cdots<
x_n$ such that for $i<n$,
the restriction of $e$ to the interval $[x_i,x_{i+1}]$ is an inclusion which
identifies $[x_i,x_{i+1}]$ with $C(e(x_i),e(x_{i+1}))\subs U_i$ for some $U_i$
in the atlas of $Y$. If $e$ is an inclusion, we speak of a {\em simple} line
path and identify it with the subset $L\subs Y$. A subset $A\subs Y$ will be
called {\em line-connected} if every pair of points
$x,y\in A$ is connected by a simple line path $L\subs A$.

\begin{prop}\pt\label{p8}
Let $Y$ be a geodesic manifold with atlas $\VVV$, and let $e\colon C\ra Y$
be an \'etale map. Then $C$ is line-connected.
\end{prop}

\pf Let $x,y\in C$ be given. By Lemma~\ref{l2}, there is a sequence
$U_1, \ldots,U_n\in\VVV_e$ with $x\in U_1, y\in U_n$, and $U_i\cap
U_{i+1}\not= \varnothing$ for $i<n$. Choose $x_i\in U_i\cap U_{i+1}$
for $i<n$. With $x_0:= x$ and $x_n:=y$, the $C(x_i,x_{i+1})$
constitute a line path $e\colon L\ra Y$ in $C$ which connects $x$
and $y$. Assume that the interval $[x,x_i]\subs L$ maps onto a
simple line path $L'$. If $C(x_i,x_{i+1})$ intersects $L'$ in a
point $\not=x_i$, there is a largest $z\in C(x_i,x_{i+1})$ with
property. Thus, if $z'$ denotes the corresponding point on $L'$, we
can replace the interval $[z',z]$ by $\{z\}$ and attach the segment
$C(z,x_{i+1})$. After less than $n$ modifications, we get a simple
line path between $x$ and $y$. \bx

By (G2), we have the following

\vspace{3mm}
\noindent {\bf Corollary.} \it Let $Y$ be a geodesic manifold. Every geodesic
with end points $x,y\in Y$ is a line path.  \rm

\vspace{3mm}
In particular, a simple geodesic with end points $x,y\in Y$ is just a minimal
connected set $C\subs Y$ with $x,y\in C$ which is locally convex.

Let $Y$ be a geodesic manifold. For $x,y\in Y$, we define a {\em simple arc}
between $x$ and $y$ to be a subspace $A\subs Y$ which is a linear continuum
with end points $x$ and $y$. We fix a
linear order on $A$ such that $x$ becomes the smallest element and denote the
set of all such $A$ by $\Omega_Y(x,y)$. In particular, every simple line path
between $x$ and $y$ belongs to $\Omega_Y(x,y)$. Clearly, every $A\in
\Omega_Y(x,y)$ admits an inscribed line path $L$ between $x$ and $y$.
Although there is no concept of length at our disposal, the intuition that
$L$ is ``shorter'' than $A$ can be expressed by the inclusion $\overline{C(L)}
\subs\overline{C(A)}$. Thus it is natural to define a preordering on
$\Omega_Y(x,y)$ by
\begin{equation}\label{15}
A\prec B :\Equ\;\overline{C(A)}\subs\overline{C(B)}.
\end{equation}
If $A\prec B$ holds for a pair $A,B\in \Omega_Y(x,y)$, we say that
$A$ is a {\em straightening} of $B$. Define $B\in\Omega_Y(x,y)$ to
be {\em minimal} if $A\prec B$ implies $B\prec A$ for all
$A\in\Omega_Y(x,y)$. We have the following straightening theorem 
which justifies the term ``geodesic'' manifold in
Definition~\ref{d7}.
\begin{thm}\pt\label{t1}
Let $Y$ be a geodesic manifold. Every simple arc $A\in\Omega_Y(x,y)$ in $Y$
can be straightened to a minimal $C\in\Omega_Y(x,y)$. A simple arc
$A\in\Omega_Y(x,y)$ is minimal if and only if $A$ is a convex simple geodesic.
\end{thm}

\pf Let $A\in\Omega_Y(x,y)$ be given. Since $C(A)$ is connected,
$\overline{C(A)}$ is connected. Proposition~\ref{p6} implies that
$\overline{C(A)}$ is convex. So the inclusion $\overline{C(A)}\hra Y$ is
\'etale. By
Proposition~\ref{p8}, there exists a simple line path $L\subs\overline{C(A)}$
between $x$ and $y$. Hence $L\prec A$. As $L$ belongs to the convex hull of a
finite set, (G1) implies that $\overline{C(L)}$ is compact. We have to
verify that $\overline{C(L)}$ contains a minimal
$C\in\Omega_Y(x,y)$. Let $\CC$ be a chain of compact convex connected sets
$C\subs\overline{C(L)}$ with $x,y\in C$. Then $D:=\bigcap\CC$ is compact and
convex, and $x,y\in D$. We show first that every open set $V$ of $Y$ with
$D\subs V$ contains some $C\in\CC$. In fact, the set $\overline{C(L)}$
is compact, and $\bigcap_{C\in\CC}(C\setm V)=\varnothing$. Hence $C\setm V=
\varnothing$ for some $C\in\CC$.

\vspc Next we show that $D$ is connected. Suppose that there is a
disjoint union $D=D_1\sqcup D_2$ with non-empty compact sets $D_1$
and $D_2$. Then we can find open sets $U_1$ and $U_2$ in $Y$ with
$D_i\subs U_i$ such that $U_1\cap U_2=\varnothing$ (see, e.g.,
\cite{Kel}, chap. V, Theorem~8). Hence $D\subs U_1\sqcup U_2$, which
yields $C\subs U_1\sqcup U_2$ for some $C\in\CC$. Since $C$ is
connected, we can assume that $C\subs U_1$. This gives $D_2\subs
U_1\cap U_2=\varnothing$, a contradiction. Thus $D$ is connected. By
Zorn's Lemma, it follows that there exists a minimal compact convex 
connected set $C$ with $x,y\in C$. Hence $C\hra Y$ is an \'etale map
generated by $\{ x,y\}$. Therefore, (G2) implies that $C$ admits no
connected proper subset $C'\subs C$ with $x,y\in C'$. By
Proposition~\ref{p8}, it follows that $C$ is a simple line path,
whence $C\in\Omega_Y(x,y)$, and $C$ is minimal.

\vspc
In particular, we have shown that if $A\in\Omega_Y(x,y)$ is minimal, then
$A$ is a convex simple geodesic between $x$ and
$y$. Conversely, if $A\in\Omega_Y(x,y)$ is a convex simple geodesic, then
$A=\overline{C(A)}$, and thus $A$ is minimal. \bx

We conclude this section with some typical examples.

\noindent {\bf Example 5.} Let $Y$ be a geodesic manifold with atlas $\VVV$,
and let $Z$ be a closed locally convex subspace. Then $Z\hra Y$ is \'etale.
Every finite set $F$ in $Z$ is contained in a compact convex set $C$ in $Y$.
Hence $C\cap Z$ is compact and convex in $Z$. Thus $Z$ satisfies (G1). As (G2)
trivially carries over to $Z$, it follows that $Z$ is a geodesic manifold.
If $Y$ is a geodesic $q$-manifold, then so is $Z$.

\noindent {\bf Example 6.} Let $Y$ be a complete locally compact
length metric space \cite{BrH,Gr}. By the Hopf-Rinow Theorem 
(\cite{BrH}, Proposition~I.3.7), the closed metric balls in $Y$ are
compact, and any two points in $Y$ are connected by a shortest path.
It is natural to assume that $Y$ admits a basis of convex open sets
where shortest paths are unique. This provides $Y$ with a local
convexity structure $\VVV$ which satisfies (G1). Note that by
\cite{BrH}, I.3.12, the map (\ref{2}) is continuous where it is
defined.

\vspc Now let $e\colon C\ra Y$ be an \'etale map generated by
$\{x,y\}\subs C$, where $C$ is compact. Similar to the case of a
covering of length metric spaces (\cite{BrH}, Proposition~I.3.25),
the length metric $d_Y$ of $Y$ can be lifted to a length metric
$d_C$ of $C$ such that $d_C(u,v)\ge d_Y(e(u),e(v))$ for all $u,v\in
C$. (If $d_C(u,v)=0$ with $u\not=v$, a neighbourhood $U\in\VVV_e$ of
$u$ cannot contain $v$. As $U$ contains a closed neighbourhood of
$u$ in $C$, we get $d_C(u,v)>0$.) Since $C$ is compact, the
Hopf-Rinow Theorem, applied to $C$, yields a shortest path $L\subs 
C$ between $x$ and $y$. Hence $C=L$, which proves (G2). By
Proposition~\ref{p4}, $Y$ is a geodesic $q$-manifold.

\noindent {\bf Example 7.} Let $Y$ be a locally convex topological vector
space. For $x,y\in Y$, we set $C(x,y):=\{ \lambda x+(1-\lambda)y\: |\:
0\le\lambda\le 1\}$ to make $Y$ into a convexity space. For a finite set
$F\subs Y$, the closed convex hull $\overline{C(F)}$ of $F$ is contained in a
finite dimensional subspace of $Y$. Hence $\overline{C(F)}$ is compact. Thus
$Y$ satisfies (G1). Let
$e\colon C\ra Y$ be an \'etale map generated by $\{x,y\}\subs C$,
where $C$ is compact. By Proposition~\ref{p8}, $e$ is generated by a
simple line path in $C$. Hence $e(C)$ is contained in a finite dimensional
subspace of $Y$. So Example~6 applies, which proves (G2). Thus $Y$ is a
geodesic manifold. Moreover, Example~4 shows that $Y$ is a geodesic
$q$-manifold if and only if $Y$ does not contain a locally convex direct sum
$\R^{(\aleph_0)}$ as a subspace.

\section{The Lokal-global-Prinzip}

With respect to convex neighbourhoods, \'etale maps have the following
disjointness property.
\begin{prop}\pt\label{p9}
Let $Y$ be a geodesic manifold with atlas $\VVV$, and let $e\colon C\ra Y$
be an \'etale map. Assume that $U,U'\in\VVV_e$.
If $e|_{U\cup U'}$ is not injective, then $U\cap U'=\varnothing$.
\end{prop}

\pf If $e|_{U\cup U'}$ is not injective, there exist $x\in U$ and
$x'\in U'$ with $e(x)=e(x')$. Suppose that there is some $z\in U\cap
U'$. Then $x\not=z$, and $U\cap U'\cap C(x,z)$ is a convex open
subset of $C(x,z)\setm\{ x\}$. Hence there is a point $t\in C(x,z)$
with $(U\setm U')\cap C(x,z)=C(x,t)$. So the homeomorphisms
$C(x,z)\cong C(e(x),e(z))\cong C(x',z)$ give rise to a point $t'\in
U'$ with $e(t)=e(t')$ and $(U'\setm U)\cap C(x',z)=C(x',t')$.
Moreover, $D:=C(t,z)\cup C(t',z)=C(t,z)\cup\{ t'\}$ since $e|_U$ is
injective. Therefore, $D$ is not a minimally connected superset of
$\{t,z\}$. On the other hand, $D$ is compact with open subsets
$C(t,z)$ and $C(t',z)$. Hence $e|_D\colon D\ra Y$ is an \'etale map
generated by $\{t,z\}$, contrary to (G2). \bx

As an immediate consequence, the fibers of an \'etale map can be separated by
pairwise disjoint neighbourhoods.

\vspace{3mm}
\noindent {\bf Corollary 1.} \it Let $Y$ be a geodesic manifold, and let
$e\colon C\ra Y$ be an \'etale map. For a given $y\in Y$,
choose a neighbourhood $U_x\in\VVV_e$ of each $x\in f^{-1}(y)$. Then the $U_x$
are pairwise disjoint.        \rm

\vspace{3mm}
\noindent {\bf Corollary 2.} \it Let $Y$ be a geodesic manifold, and let
$e\colon C\ra Y$ be an \'etale map. Then $C$ is a Hausdorff space. \rm

\vspace{3mm}
\pf Let $x,x'\in C$
be given. If $e(x)\not= e(x')$, there are disjoint neighbourhoods of $e(x)$
and $e(x')$, and their inverse images give disjoint neighbourhoods of $x$ and
$x'$. So we can assume that $e(x)=e(x')$. Choose $U,U'\in\VVV_e$ with $x\in U$
and $x'\in U'$. By Proposition~\ref{p9}, $U\cap U'=\varnothing$. Thus $C$ is
Hausdorff. \bx

If the geodesic manifold is regular, the fibers are even discrete, which leads
to the following finiteness result.

\begin{prop}\pt\label{p10}
Let $e\colon C\ra Y$ be an \'etale map into a geodesic $q$-manifold $Y$. Then
the fibers of $e$ are finite.
\end{prop}

\pf Let $\VVV$ denote the atlas of $Y$, and let $y\in Y$ be given. For each
$x\in e^{-1}(y)$, we choose a neighbourhood $U_x\in\VVV_e$ such that the
images $e(U_x)$ are contained in a fixed $V'\in\VVV$. By the
Corollary~1, these neighbourhoods are pairwise disjoint. Without loss of
generality, we can assume that $|C|>1$. Since $C$ is a connected Hausdorff
space by Corollary~2, this implies that $C$ has no isolated points. As $e$ is
closed, the complement of $\bigcup\{ U_x\: |\: x\in e^{-1}(y)\}$ is mapped to
a closed set $A\subs Y$ with $y\notin A$. So there exists an open
neighbourhood $W\subs V'$ of $y$ with $e^{-1}(W)\subs\bigcup\{ U_x\: |\: x\in
e^{-1}(y)\}$. By the regularity of $Y$, we find a convex open neighbourhood
$V$ of $y$ with $\overline{V}\subs W$.

\vspc
For any $x\in e^{-1}(y)$, the set $U_x\cap e^{-1}(V)$ is an open neighbourhood
of $x$, hence not a singleton. Therefore, the $V_x:=e(U_x\cap e^{-1}(V))$
are convex subsets of $V$ with $|V_x|>1$ and $y\in V_x$. Choose arbitrary
$z_x\in U_x\cap e^{-1}(V)$ with $y_x:=e(z_x)\not=y$ for all $x\in e^{-1}(y)$.
Now let $Z\subs\bigcup\{C(x,z_x)\:|\:x\in e^{-1}(y)\}$ be such that
$Z\cap C(x,z_x)$ is closed in $U_x\cap e^{-1}(V)$ for every $x\in e^{-1}(y)$.
We claim that $Z$ is closed. Thus let $z\in\overline{Z}$ be given. Then
$e(z)\in\overline{e(Z)}\subs\overline{V}\subs W$. Hence
$z\in e^{-1}(W)\subs\bigcup\{U_x\: |\:\mbox{$x\in e^{-1}(y)$}\}$, which
yields $z\in Z$. Thus $Z$ is closed. Since $e$ is closed, this implies
that $S(y):=\bigcup\{ C(y,y_x)\: |\: x\in e^{-1}(y)\}$ is closed and carries
the finest topology such that the maps $C(y,y_x)\hra S(y)$ are continuous for
all $x\in e^{-1}(y)$.

\vspc Suppose that $e^{-1}(y)$ is infinite. By Ramsey's Theorem 
\cite{Ram}, there must be an infinite subset $E$ of $e^{-1}(y)$ such
that either $C(y,y_u)\cap C(y,y_v)=\{y\}$ for all pairs of different
$u,v\in E$, or $C(y,y_u)\cap C(y,y_v)\not=\{y\}$ for different
$u,v\in E$. The first case is impossible since $V$ is star-finite by
Definition~\ref{d7}. Otherwise, there is a point $y'\in V\setm\{y\}$
and a set $Z\subs\bigcup\{C(x,z_x)\:|\:x\in e^{-1}(y)\}$ with
$|Z\cap C(x,z_x)|=1$ for all $x\in E$ such that $e(Z)$ is an
infinite non-closed subset of $C(y,y')$. Since $Z$ is closed, this
gives a contradiction. \bx

As a consequence, the geodesic structure of a geodesic $q$-manifold can be
lifted along \'etale maps.

\begin{thm}\pt\label{t2}
Let $e\colon C\ra Y$ be an \'etale map into a geodesic $q$-manifold $Y$ with
atlas $\VVV$. Then $C$ is a geodesic $q$-manifold with atlas $\VVV_e$.
\end{thm}

\pf By Corollary~2 of Proposition~\ref{p9}, $C$ is a Hausdorff space. We show
first that $C$ is regular. Let $U_x\in\VVV_e$ be a neighbourhood
of $x\in C$. We choose neighbourhoods $U_z\in\VVV_e$ for all $z$ in the fiber
of $y:=e(x)$. By Corollary~1 of Proposition~\ref{p9}, the $U_z$ are pairwise
disjoint. Since $Y$ is regular and $e$ closed, there is a closed
neighbourhood $V$ of $y$ with $e^{-1}(V)\subs\bigcup\{ U_z\: |\: z\in e^{-1}(y)
\}$. Hence
$$U_x\cap e^{-1}(V)=e^{-1}(V)\setm\bigcup\bigl\{ U_z\: |\: z\in
e^{-1}(y)\setm\{ x\}\bigr\}$$
is a closed neighbourhood of $x$. Thus $C$ is regular.

\vspc
Let $F\subs C$ be finite. Then $\overline{C(e(F))}$ is compact. By
Proposition~\ref{p10}, the fibers of $e$ are compact. Hence
$e^{-1}(\overline{C(e(F))})$ is compact by \cite{Bou}, chap. I.10,
Proposition~6. Furthermore, $e^{-1}(\overline{C(e(F))})$ is convex with respect
to $\VVV_e$. Therefore, the closed subset $\overline{C(F)}$ of
$e^{-1}(\overline{C(e(F))})$ is compact. This proves (G1) for $C$.

\vspc
Next let $e'\colon C'\ra C$ be an \'etale map with $C'$ compact which is
generated by $\{ x,y\}\subs C'$. Then $ee'$ is \'etale and generated by
$\{ x,y\}$. Hence $C'$ is minimal among the connected sets $B\subs C'$ with
$x,y\in B$. Thus $C$ satisfies (G2).

\vspc
Finally, let $S(x,E):=\bigcup\{ C(x,z)\:|\:z\in E\}$ be a closed star in some
$U\in\VVV_e$. Since $C$ is regular, we find a closed convex neighbourhood
$U'\subs U$ of $x$. By Proposition~\ref{p3}, this implies that $S(x,E)\cap U'$
is a star in $U$ which is closed in $C$.
Therefore, $e(S(x,E)\cap U')$ is a closed star in some $V\in\VVV$. So $E$ is
finite, which proves that $C$ is a geodesic $q$-manifold. \bx

Now we are ready to prove our main result which essentially states that the
image of an \'etale map is weakly convex in the following sense
(cf. \cite{BOR1}, Definition~2.16).
\begin{Definition}\pt\label{d9}
\rm Let $Y$ be a geodesic manifold. We call a subset $A\subs Y$ {\em weakly
convex} if every pair of points $x,y\in A$ can be connected by a geodesic.
\end{Definition}

The following theorem extends previous versions of the
Lokal-global-Prinzip for convexity of maps (see 
\cite{CDM,HNP,BOR1,BOR2}).

\begin{thm}\pt\label{t3}
Let $f\colon X\ra Y$ be a locally convex continuous map from a connected
topological
space $X$ to a geodesic $q$-manifold $Y$. Assume that $f^{\#}$ is closed.
Then $f(X)$ is weakly convex.
\end{thm}

\pf Let $\VVV$ be the atlas of $Y$. By Proposition~\ref{p7}, the map $f^{\#}$
again is locally convex, and Proposition~\ref{p5} implies that $f^{\#}$ is
\'etale. By Theorem~\ref{t2}, it follows that $X^f$ is a geodesic manifold.
For $z,z'\in X^f$, Proposition~\ref{p8} shows that there is a connecting
simple line path $L$ between $z$ and $z'$. Theorem~\ref{t1} shows that $L$ can
be straightened to a convex simple geodesic $C$. Thus $f^{\#}|_C\colon C\ra Y$
is a geodesic between $f^{\#}(z)$ and $f^{\#}(z')$. Hence $f(X)$ is
weakly convex. \bx

In the special case where $f$ is an inclusion $X\hra Y$, the preceding proof
yields

\vspace{3mm}
\noindent {\bf Corollary.} \it Let $C$ be a closed connected locally convex
subset of a geodesic manifold $Y$. Then $C$ is weakly convex.       \rm

\vspace{3mm}
\pf By Example~5, $C$ is a geodesic manifold, and $C\hra Y$ is \'etale. As in
the proof of Theorem~\ref{t3}, this implies that $C$ is weakly convex. \bx

\noindent {\bf Remarks. 1.} If $f$ is closed, then $f^{\#}$ is closed. However,
the latter condition is much weaker. For example, if $f$ is a local
homeomorphism, then $f^{\#}$ is identical, but $f$ need not be closed.

{\bf 2.} The preceding corollary extends Klee's generalization of a 
classical result due to Tietze \cite{Ti} and Nakajima (Matsumura)
\cite{Na}. Klee's Theorem \cite{Klee} states that the above 
corollary holds in a locally convex topological vector space $Y$.
Note that the usual proof of Klee's Theorem rests on the linear 
structure of $Y$, while the corollary of Theorem~\ref{t3} merely
depends on a local convexity structure in the sense of
Definition~\ref{d4}.

\vspace{5mm} \noindent {\Large\bf Acknowledgements}

This work was started on the occasion of Tudor S. Ratiu's sixtieth
birthday. The second author would like to thank the late Jerrold E.
Marsden (without his promise the work would never have been born!),
Tudor S. Ratiu, and Karl-Hermann Neeb for helpful conversations.
Special thanks are extended to Erwin-Schrödinger International
Institute for Mathematical Physics in Vienna for stimulating
atmosphere, Ralph L. Cohen for his invitation to attend the
Copenhagen Topology Conference (January 8-10, 2010) in Copenhagen,
Denmark, that has opened new possibilities, Alan D. Weinstein for
support, and Dennis P. Sullivan for awakening interest and
stimulating conversations especially in Topology.

Parts of this work were presented at the Young Topologists Meeting
2010 in K{\o}benhavn, Denmark, June 16-20, 2010, at Geometry,
Mechanics, and Dynamics: A workshop celebrating the 60th birthday of
Tudor Ratiu at CIRM Luminy, France, July 12-16, 2010, and at
Geometric and Algebraic Structures in Mathematics: A conference to
celebrate Dennis Sullivan's 70th birthday, held at Simons Center for
Geometry and Physics, Stony Brook, USA, May 26 - June 4, 2011. Kind
invitations, hospitalities, and grants by the organizers of the
events, Stony Brook University, Simons Foundation, Aarhus and
K{\o}benhavns Universitet, are highly acknowledged.

---

\end{document}